\documentstyle[amssymb]{article}
\begin{document}

\noindent
\centerline{\bf THE NONSTANDARD DEFORMATION U$'_{\bf q}$(so$_{\bf n}$)}
\centerline{\bf FOR q A ROOT OF UNITY}
\vskip 10 pt

\centerline{N. Z. IORGOV AND A. U. KLIMYK
\footnote{
The research described in this publication was made possible
in part by CRDF Grant UP1-309 and by DFFD Grant 1.4/206.}}

\vskip  25 pt

\begin{abstract}
We describe properties of the nonstandard $q$-deformation
$U'_q({\rm so}_n)$ of the universal
enveloping algebra $U({\rm so}_n)$ of the Lie algebra ${\rm so}_n$
which does not coincide with the Drinfeld--Jimbo quantum algebra
$U_q({\rm so}_n)$. In particular, it is shown that there exists
an isomorphism from $U'_q({\rm so}_n)$ to $U_q({\rm sl}_n)$ and that
finite dimensional irreducible representations of $U'_q({\rm so}_n)$
separate elements of this algebra.
Irreducible representations of the algebras $U'_q({\rm so}_n)$
for $q$ a root of unity $q^p=1$ are given. The main class of
these representations act
on $p^N$-dimensional linear space (where $N$ is a number of positive
roots of the Lie algebra ${\rm so}_n$) and are given by
$r={\rm dim}\, {\rm so}_n$ complex parameters. Some classes of degenerate
irreducible representations are also described.
\end{abstract}

\vskip 35 pt

\centerline{\sc 1. Introduction}
\medskip

\noindent
Quantum orthogonal groups, quantum Lorentz groups and their
corresponding quantum algebras are of special interest for modern
mathematical physics [1-3].
M. Jimbo [4] and V. Drinfeld [5] defined $q$-deformations
(quantum algebras) $U_q(g)$ for all simple complex Lie algebras $g$
by means of Cartan subalgebras and root subspaces (see also [6] and [7]).
Reshetikhin,
Takhtajan and Faddeev [8] defined quantum algebras $U_q(g)$ in terms
of the quantum $R$-matrix satisfying the quantum Yang--Baxter
equation. However, these approaches do not
give a satisfactory presentation of the quantum algebra
$U_q({\rm so}(n,{\bf C}))$ from a viewpoint of some problems in quantum
physics and representation theory.
When considering representations of the quantum
groups $SO_q(n+1)$ and $SO_q(n,1)$
we are interested in reducing them onto the quantum subgroup $SO_q(n)$.
This reduction would give an analogue of the Gel'fand--Tsetlin basis for
these representations. However, definitions of quantum algebras
mentioned above do not allow the inclusions
$U_q({\rm so}(n+1, {\bf C}))\supset U_q({\rm so}(n,{\bf C}))$
and $U_q({\rm so}_{n,1})\supset U_q({\rm so}_n)$. To be able to exploit such
reductions we have to consider $q$-deformations of the Lie algebra
${\rm so}(n+1,{\bf C})$ defined in terms of the generators
$I_{k,k-1}=E_{k,k-1}-E_{k-1,k}$  (where $E_{is}$ is the matrix with
elements $(E_{is})_{rt}=\delta _{ir} \delta _{st})$ rather than
by means of Cartan subalgebras and root elements.
To construct such deformations we have to deform trilinear relations
for elements $I_{k,k-1}$ instead of Serre's relations (used in the case of
Jimbo's quantum algebras). As a result, we obtain the associative algebra
which will be denoted as $U'_q({\rm so}(n,{\bf C})).$

     These $q$-deformations were first constructed  in [9].
They permit one to construct the reductions of $U'_q({\rm so}_{n,1})$ and
$U'_q({\rm so}_{n+1})$ onto $U'_q({\rm so}_n)$.
The $q$-deformed algebra $U'_q({\rm so}(n,{\bf C}))$ leads for $n=3$ to
the $q$-deformed algebra $U'_q({\rm so}(3,{\bf C}))$ defined by
D. Fairlie [10]. The cyclically symmetric algebra, similar to Fairlie's
one, was also considered somewhat earlier by Odesskii [11].
The algebra $U'_q({\rm so}(3,{\bf C}))$ allows us to construct the
noncompact quantum algebra $U'_q({\rm so}_{2,1})$. The algebra
$U'_q({\rm so}(4,{\bf C}))$ is a $q$-deformation of the Lie
algebra ${\rm so}(4,{\bf C})$ given by means of usual bilinear
commutation relations between the elements $I_{ji}$, $1\le i<j\le 4$.
In the case of the classical
Lie algebra ${\rm so}(4,{\bf C})$ one has
${\rm so}(4,{\bf C})={\rm so}(3,{\bf C})+{\rm so}(3,{\bf C})$,
while in the case of our $q$-deformation $U'_q({\rm so}(4,{\bf C}))$
this is not the case.

In the classical case, the imbedding $SO(n)\subset SU(n)$
(and its infinitesimal analogue) is of great importance for nuclear
physics and in the theory of Riemannian symmetric
spaces. It is well known
that in the framework of Drinfeld--Jimbo quantum groups and algebras
one cannot construct the corresponding embedding. The algebra
$U'_q({\rm so}(n,{\bf C}))$ allows to define such an embedding [12],
that is, it is possible to define the embedding
$U'_q({\rm so}(n,{\bf C}))\subset U_q({\rm sl}_n)$,
where $U_q({\rm sl}_n)$ is the Drinfeld-Jimbo quantum algebra.

As a disadvantage of the algebra $U'_q({\rm so}(n,{\bf C}))$ we have
to mention the difficulties with Hopf algebra structure. Nevertheless,
$U'_q({\rm so}(n,{\bf C}))$ turns out to be a coideal in
$U_q({\rm sl}_n)$ (see [12]) and this fact allows us to consider tensor
products of finite dimensional irreducible representations of
$U'_q({\rm so}(n,{\bf C}))$ for many interesting cases.

For convenience, below we denote the Lie algebra ${\rm so}(n,{\bf C})$ by
${\rm so}_n$ and the $q$-deformed algebra $U'_q({\rm so}(n,{\bf C}))$ by
$U'_q({\rm so}_n)$.

Finite dimensional irreducible representations of the algebra
$U'_q({\rm so}_n)$ were constructed in [9]. The formulas
of action of the generators of $U'_q({\rm so}_n)$ upon the
basis (which is a $q$-analogue of the Gel'fand--Tsetlin basis) are
given there. A proof of these formulas and some their corrections were
given in [13]. However,
finite dimensional irreducible representations described in [9] and [13]
are representations of the classical type. They are $q$-deformations of the
corresponding irreducible representations of the Lie algebra
${\rm so}_n$, that is, at $q\to 1$ they turn into representations
of ${\rm so}_n$.

The algebra $U'_q({\rm so}_n)$ has other classes of finite
dimensional irreducible representations which have no classical analogue.
These representations are singular at the limit $q\to 1$.
They are described in [14].
Note that the description of these
representations for the algebra $U'_q({\rm so}_3)$ is given in
[15]. A classification of irreducible $*$-representations of real forms
of the algebra $U'_q({\rm so}_3)$ is given in [16].

The aim of this paper is to give irreducible
representations of the algebra $U'_q({\rm so}_n)$ in the case when $q$
is a root of unity. We prove that in this case all irreducible
representations of $U'_q({\rm so}_n)$ are finite dimensional. In
order to prove the corresponding theorem we prove an analogue of
the Poincar\'e--Birkhoff--Witt theorem for $U'_q({\rm so}_n)$ (this analogue
was anounced in [17]) and use central elements of this algebra
for $q$ a root of unity (they are derived in [18]).

For construction of irreducible representations of $U'_q({\rm so}_n)$
for $q$ a root of unity, we use the method of D. Arnaudon and A.
Chakrabarti [19] for construction of
irreducible representations of the quantum algebra $U_q({\rm sl}_n)$
when $q$ is a root of unity. If $q^p=1$ and $p$ is an odd integer, then
we construct the series of
irreducible representations of $U'_q({\rm so}_n)$ which act on
$p^N$-dimensional vector space (where $N$ is the number of positive roots
of the Lie algebra ${\rm so}_n$) and are given by $r={\rm dim}\, {\rm so}_n$
complex parameters. These representations are irreducible for generic
values of these parameters. These representations constitute the main
class of irreducible representations of $U'_q({\rm so}_n)$.
For some special values of the representation parameters in ${\Bbb C}^r$
the representations are reducible. These reducible representations give
many other classes of (degenerate)
irreducible representations which are given
by less number of parameters or by parameters, values of which cover
subsets of ${\Bbb C}^r$ of Lebesgue measure 0. As in the case of
irreducible representations of the quantum algebra $U_q({\rm sl}_n)$,
it is difficult to enumerate all irreducible representations of these
classes. However, we give some most important classes of these
degenerate representations. In particular, we give $2^{n-1}$ classes
of these representations, which are an analogue of the nonclassical type
irreducible representations of $U'_q({\rm so}_n)$ for $q$ not a root
of unity.
\bigskip

\centerline{\sc 2. The $q$-deformed algebra $U'_q({\rm so}_n)$}
\medskip

\noindent
The Drinfeld--Jimbo algebra $U_q({\rm so}_n)$ is obtained by deforming
Serre's relations for generating elements $E_1,\cdots , E_l$,
$F_1,\cdots , F_l$, $H_1,\cdots , H_l$ of $U({\rm so}_n)$ (see [20]).
In order
to obtain $U'_q({\rm so}_n)$ we have to take determining relations for
the generating elements $I_{21}$, $I_{32},\cdots ,I_{n,n-1}$ of
$U({\rm so}_n)$ (they do not coincide with $E_j$, $F_j$, $H_j$) and to
deform these relations. The elements $I_{21}$, $I_{32},\cdots ,I_{n,n-1}$
belong to the basis $I_{ij}$, $i>j$, of the Lie algebra ${\rm so}_n$.
The matrices $I_{ij}$, $i>j$, are defined as $I_{ij}=E_{ij}-E_{ji}$,
where $E_{ij}$ is the matrix with entries $(E_{ij})_{rs}=
\delta _{ir}\delta _{js}$. The universal enveloping algebra
$U({\rm so}_n)$ is generated by a part of the basis elements $I_{ij}$,
$i>j$, namely, by the elements $I_{21}$, $I_{32},\cdots ,I_{n,n-1}$.
These elements satisfy the relations
$$
I^2_{i,i-1}I_{i+1,i}-2I_{i,i-1}I_{i+1,i}I_{i,i-1} +
I_{i+1,i}I^2_{i,i-1} =-I_{i+1,i}, $$
$$
I_{i,i-1}I^2_{i+1,i}-2I_{i+1,i}I_{i,i-1}I_{i+1,i} +
I^2_{i+1,i}I_{i,i-1} =-I_{i,i-1}, $$
$$
I_{i,i-1}I_{j,j-1}- I_{j,j-1}I_{i,i-1}=0\ \ \ \ {\rm for}\ \ \ \
|i-j|>1.
$$
The following theorem is true [21] for the
universal enveloping algebra $U({\rm so}_n)$.
\medskip

\noindent
{\bf Theorem 1.} {\it The universal enveloping algebra $U({\rm so}_n)$
is isomorphic to the complex associative algebra (with a unit element)
generated by the elements $I_{21}$, $I_{32},\cdots ,I_{n,n-1}$
satisfying the above relations.}
\medskip

We make the $q$-deformation of these relations by
$2\to [2]:=(q^2-q^{-2})/(q-q^{-1})=q+q^{-1}$.
As a result, we obtain the complex associative algebra
generated by elements $I_{21}$, $I_{32},\cdots ,I_{n,n-1}$ satisfying
the relations
$$
I^2_{i,i-1}I_{i+1,i}-(q+q^{-1})I_{i,i-1}I_{i+1,i}I_{i,i-1} +
I_{i+1,i}I^2_{i,i-1} =-I_{i+1,i}, \eqno (1) $$
$$
I_{i,i-1}I^2_{i+1,i}-(q+q^{-1})I_{i+1,i}I_{i,i-1}I_{i+1,i} +
I^2_{i+1,i}I_{i,i-1} =-I_{i,i-1}, \eqno (2) $$
$$
I_{i,i-1}I_{j,j-1}- I_{j,j-1}I_{i,i-1}=0\ \ \ \ {\rm for}\ \ \ \
|i-j|>1.  \eqno (3)
$$
This algebra was introduced by us in [9] and is denoted by
$U'_q({\rm so}_n)$.

We wish to formulate and to prove for the algebra $U'_q({\rm so}_n)$
an analogue of the Poincar\'e--Birkhoff--Witt theorem.
For this we determine (see [22] and [23])
in $U'_q({\rm so}_n)$ elements analogous to the
matrices $I_{ij}$, $i>j$, of the Lie algebra ${\rm so}_n$. In order
to give them we use the
notation $I_{k,k-1}\equiv I^+_{k,k-1}\equiv I^-_{k,k-1}$.
Then for $k>l+1$ we define recursively
$$
I^+_{kl}:= [I_{l+1,l},I_{k,l+1}]_{q}\equiv
q^{ 1/2}I_{l+1,l}I_{k,l+1}-
q^{- 1/2}I_{k,l+1}I_{l+1,l}, \eqno (4)
$$
$$
I^-_{kl}:= [I_{l+1,l},I_{k,l+1}]_{q^{-1}}\equiv
q^{-1/2}I_{l+1,l}I_{k,l+1}-
q^{1/2}I_{k,l+1}I_{l+1,l}.
$$
The elements $I^+_{kl}$, $k>l$, satisfy the commutation relations
$$
[I^+_{ln},I^+_{kl}]_q=I^+_{kn},\ \
[I^+_{kl},I^+_{kn}]_q=I^+_{ln},\ \
[I^+_{kn},I^+_{ln}]_q=I^+_{kl} \ \ \
{\rm for}\ \ \  k>l>n, \eqno (5) $$
$$
[I^+_{kl},I^+_{nr}]=0\ \ \ \ {\rm for}\ \ \
k>l>n>r\ \ {\rm and}\ \ k>n>r>l, \eqno (6)
$$
$$
[I^+_{kl},I^+_{nr}]_q=(q-q^{-1})
(I^+_{lr}I^+_{kn}-I^+_{kr}I^+_{nl}) \ \ \ {\rm for}\ \ \
k>n>l>r.   \eqno(7)
$$
For $I^-_{kl}$, $k>l$, the commutation relations are obtained
from these relations by replacing
$I^+_{kl}$ by $I^-_{kl}$ and $q$ by $q^{-1}$.

The algebra $U'_q({\rm so}_n)$ can be considered as an
associative algebra (with unit element) generated by $I^+_{kl}$,
$1\le l<k\le n$, satisfying the relations (5)--(7). Really, using
the relations (4) we can reduce the relations (5)--(7) to the
relations (1)--(3) for $I_{21}$, $I_{32},\cdots ,I_{n,n-1}$
(for the case of the algebra $U'_q({\rm so}_3)$ this reduction is
simple and is given, for example, in [15]). Similarly,
$U'_q({\rm so}_n)$ is an associative algebra generated by
$I^-_{kl}$, $1\le l<k\le n$, satisfying the corresponding relations.

Now the Poincar\'e--Birkhoff--Witt theorem for the
algebra $U'_q({\rm so}_n)$ can be formulated as follows.
\medskip

\noindent
{\bf Theorem 2.} {\it The elements
$$
{I_{21}^+}^{m_{21}}{I_{31}^+}^{m_{31}}\cdots {I_{n1}^+}^{m_{n1}}
{I_{32}^+}^{m_{32}} {I_{42}^+}^{m_{42}} \cdots {I_{n2}^+}^{m_{n2}}
\cdots {I_{n,n-1}^+}^{m_{n,n-1}},\ \ \ \  m_{ij}=0,1,2, \cdots ,
$$
form a basis of the algebra $U'_q({\rm so}_n)$.
This assertion is true if $I^+_{ij}$ are replaced by the
corresponding elements $I^-_{ij}$.}
\medskip

\noindent
{\sl Proof.}
The proof of this theorem is essentially based on
Bergman's Diamond Lemma [24]. All the terms which we use in this
proof without explanation are defined there. Let us consider
$U'_q({\rm so}_n)$ as an associative algebra, generated by the
elements $I^+_{kl}$, $1\le l < k \le n$, satisfying the relations
(5)--(7). These relations can be presented in the form of
reduction rules \setcounter{equation}{7}
\begin{eqnarray}
I^+_{kl} I^+_{ml} & = & q^{-1}  I^+_{ml}  I^+_{kl} + q^{-1/2}
I^+_{km}, \qquad\qquad \qquad\qquad\qquad k>m>l\ ,\\ I^+_{km}
I^+_{ml} & = & q  I^+_{ml}  I^+_{km} - q^{1/2} I^+_{kl},
\qquad\qquad\qquad\qquad\qquad\qquad k>m>l\ ,\\ I^+_{km} I^+_{kl}
& = & q^{-1}  I^+_{kl}  I^+_{km} + q^{-1/2} I^+_{ml},
\qquad\qquad\qquad\qquad\qquad k>m>l\ ,\\ I^+_{kl} I^+_{mp} & = &
I^+_{mp}I^+_{kl},  \qquad \qquad\qquad\qquad\qquad\qquad
\qquad\qquad k>l>m>p\ ,\\ I^+_{kl} I^+_{mp} & = &
I^+_{mp}I^+_{kl},  \qquad\qquad\qquad\qquad\qquad
\qquad\qquad\qquad m>k>l>p\ ,\\ I^+_{kl} I^+_{mp} & = &
I^+_{mp}I^+_{kl} + (q-q^{-1}) (I^+_{lp} I^+_{km} - I^+_{kp}
I^+_{ml})\ ,  \quad\qquad k>m>l>p\ .
\end{eqnarray}
Every element of the algebra $U'_q({\rm so}_n)$ can be presented as a linear
combination of monomials of the noncommuting elements $I^+_{kl}$,
$l<k$. If some monomial contains as a submonomial the left-hand
side of some of formulas (8)--(13), then this submonomial must
be replaced by the corresponding right-hand side.

In order to show that the described procedure of reductions will
terminate, we introduce the total ordering in the set of all
monomials. We set $I_{m,p}^+\prec I_{k,l}^+$ if either $p<l$ or
both $p=l$ and $m<k$. Then we say $X\prec Y$ if the length (the
number of generators) of monomial $X$ is less than the length of
monomial $Y$ or if their lengths are equal, but $X$ is less than
$Y$ in the sense of lexicographical ordering with respect to the
ordering of generators. The introduced ordering has an obvious
property: condition $X\prec Y$ implies $A X B\prec A Y B$ for
arbitrary two monomials $A$ and $B$. Since the left-hand side of
any rule from the reduction system (8)--(13) is greater than any
monomial in the corresponding right-hand side, the procedure of
reductions must terminate. The basis monomials of the statement of
the theorem are exactly that monomials which can not be reduced
more.

We need only to show that the result of reductions does not depend
on the order of using reduction rules. The Diamond Lemma claims
that this requirement will be fulfilled if one shows that all
ambiguities that arise in the reduction system (8)--(13) are
resolvable. It is easy to see, that all these ambiguities are {\it
overlap} ambiguities and appear when one considers monomials such
as $I^+_{i_1,i_2} I^+_{i_3,i_4} I^+_{i_5,i_6}$, where
$I^+_{i_1,i_2} \succ I^+_{i_3,i_4} \succ I^+_{i_5,i_6}$. The proof
of the resolvability of arising ambiguities is the same (up to
replacement of indices) for all the initial monomials having the same
ordering of indices. In this case, we say that ambiguities
are of the same type. It is enough to prove resolvability only for
one representative from the set of all monomials with some fixed type of
ambiguity.

Let us demonstrate resolvability of some concrete type of
ambiguity. Consider the reduction of the monomials $I^+_{i_1,i_2}
I^+_{i_3,i_4} I^+_{i_5,i_6}$ with the following ordering of
indices: $i_1>i_2=i_3=i_5>i_4>i_6$. Choose the representative
$I^+_{43} I^+_{32} I^+_{31}$ from this set of monomials. It can be
reduced in two different ways (over signs ``='' we write down
reduction rules from the reduction system (8)--(13) which must be used):
\[
(I^+_{43} I^+_{32}) I^+_{31} \stackrel{(9)}{=} (q I^+_{32}
I^+_{43} - q^{1/2} I^+_{42}) I^+_{31} \stackrel{(9)}{=} q I^+_{32}
(q I^+_{31} I^+_{43} - q^{1/2} I^+_{41})- q^{1/2} I^+_{42}
I^+_{31}=
\]
\[
\qquad\quad\stackrel{(12),(13)}{=} q^2 I^+_{32} I^+_{31}I^+_{43} -
q^{3/2} I^+_{41}I^+_{32} - q^{1/2} \Bigl(I^+_{32}
I^+_{42}+(q-q^{-1}) (I^+_{21} I^+_{43} -I^+_{41} I^+_{32})\Bigr)
=
\]
\[
\qquad\quad\stackrel{(10)}{=} q I^+_{31} I^+_{32} I^+_{43} -
q^{-1/2} I^+_{41} I^+_{32}  - q^{1/2} I^+_{31} I^+_{42} +
q^{-1/2} I^+_{21} I^+_{43}\ ,
\]
\[
I^+_{43} (I^+_{32} I^+_{31}) \stackrel{(10)}{=} I^+_{43} (q^{-1}
I^+_{31}I^+_{32} +
 q^{-1/2} I^+_{21}) \stackrel{(9),(11)}{=}
q^{-1} (q I^+_{31} I^+_{43} - q^{1/2} I^+_{41})I^+_{32}+ \qquad
\]
\[
\qquad\quad+q^{-1/2} I^+_{21} I^+_{43}\stackrel{(9)}{=} I^+_{31}
(q I^+_{32} I^+_{43} - q^{1/2} I^+_{42}) - q^{-1/2}
I^+_{41}I^+_{32} + q^{-1/2} I^+_{21} I^+_{43}\ .
\]
The results of these two reductions
are the same. This means that arising {\it overlap}
ambiguity is resolvable.

It is easy to show by exactly the same calculation, that all the other 65
types of ambiguities which arise in the reduction system (8)--(13)
are also resolvable. Therefore, all the conditions of the Diamond Lemmma
are fulfilled and theorem is proved.
\bigskip

\centerline{\sc 3. The isomorphism $U'_q({\rm so}_n)\to U_q({\rm sl}_n)$}
\medskip

\noindent
The algebra $U'_q({\rm so}_n)$ can be embedded into the Drinfeld--Jimbo
quantum algebra $U_q({\rm sl}_n)$ (see [12]). This quantum algebra is
generated by
the elements $E_i$, $F_i$, $K_i^{\pm 1}=q^{\pm H_i}$,
$i=1,2,\cdots ,n-1$, satisfying the relations
$$
K_iK_j=K_jK_i,\ \ \ K_iK_i^{-1}=K_i^{-1}K_i=1,
$$  $$
K_iE_jK_i^{-1}=q^{a_{ij}}E_j,\ \ \
K_iF_jK_i^{-1}=q^{-a_{ij}}F_j,
$$  $$
[E_i,F_j]=\delta _{ij}\frac{K_i-K_i^{-1}}{q-q^{-1}},
$$  $$
E^2_iE_{i\pm 1}-(q+q^{-1})E_iE_{i\pm 1}E_i+E_{i\pm 1}E^2_i=0,
$$   $$
F^2_iF_{i\pm 1}-(q+q^{-1})F_iF_{i\pm 1}F_i+F_{i\pm 1}F^2_i=0,
$$   $$
[E_i,E_j]=0,\ \ [F_i,F_j]=0 \ \ \ {\rm for}\ \ \ |i-j|>1,
$$
where $a_{ij}$ are elements of the Cartan matrix of the Lie algebra
${\rm sl}_n$.

In order to prove Theorem 3 below, we note that there exists a
one-to-one correspondence between the basis elements of the algebra
$U'_q({\rm so}_n)$ from Theorem 2 and the basis elements of the subalgebra
${\frak N}^-$ of the quantum algebra $U_q({\rm sl}_n)$, generated by
$F_1,F_2,\cdots ,F_{n-1}$. The last basis elements are constructed by
means of the following ordering of positive roots of the Lie algebra
${\rm sl}_n$:
$$
\beta _{12},\beta _{13},\cdots ,\beta _{1n}, \beta _{23},\cdots ,
\beta _{2n},\cdots , \beta _{n-2,n-1}, \beta _{n-2,n}, \beta _{n-1,n},
\eqno (14)
$$
where $\beta _{ij}=\alpha _i+\alpha _{i+1}+\cdots +\alpha _{j-1}$ and
$\alpha _k$ are simple roots. (This ordering is the same as in Theorem 2.)
To every root of this set there corresponds
the element $F_{\beta _{ij}}\in {\frak N}^-$ (see, for example, [7]).
Then according to the Poincar\'e--Birkhoff--Witt theorem for the algebra
${\frak N}^-$ (see, [7], subsection 6.2.3) the elements
$$
F^{m_{12}}_{\beta _{12}} F^{m_{13}}_{\beta _{13}}\cdots
F^{m_{1n}}_{\beta _{1n}}\cdots F^{m_{n-1,n}}_{\beta _{n-1,n}},\ \ \ \
m_{ij}=0,1,2,\cdots , \eqno (15)
$$
(the order of $\beta _{ij}$ is the same as in (14)) form a basis of
${\frak N}^-$. Then the mapping
$$
F^{m_{12}}_{\beta _{12}} F^{m_{13}}_{\beta _{13}}\cdots
F^{m_{n-1,n}}_{\beta _{n-1,n}} \to
{I^-_{21}}^{m_{12}} {I^-_{31}}^{m_{13}}\cdots
{I^-_{n,n-1}}^{m_{n-1,n}}  \eqno (16)
$$
is the one-to-one correspondence between basis elements in ${\frak N}^-$
and in $U'_q({\rm so}_n)$ which will be denoted by ${\cal T}$.

Similarly, to every root $\beta _{ij}$ from (14) there corresponds
the element $E_{\beta _{ij}}$ of the subalgebra
${\frak N}^+\subset U_q({\rm sl}_n)$, generated by
$E_1,E_2,\cdots ,E_{n-1}$. The elements
$$
E^{m_{12}}_{\beta _{12}} E^{m_{13}}_{\beta _{13}}\cdots
E^{m_{1n}}_{\beta _{1n}}\cdots E^{m_{n-1,n}}_{\beta _{n-1,n}},\ \ \ \
m_{ij}=0,1,2,\cdots ,
$$
(the order of $\beta _{ij}$ is the same as in (14)) form a basis of
${\frak N}^+$.

The formulas
$$
{\rm deg}\, (F^{m_{12}}_{\beta _{12}} F^{m_{13}}_{\beta _{13}}\cdots
F^{m_{n-1,n}}_{\beta _{n-1,n}})=-(m_{12}\beta _{12}+
m_{13}\beta _{13}+\cdots m_{n-1,n}\beta _{n-1,n}),
$$  $$
{\rm deg}\, (E^{m_{12}}_{\beta _{12}} E^{m_{13}}_{\beta _{13}}\cdots
E^{m_{n-1,n}}_{\beta _{n-1,n}})=m_{12}\beta _{12}+
m_{13}\beta _{13}+\cdots m_{n-1,n}\beta _{n-1,n},
$$   $$
{\rm deg}\, (H_1^{m_1}\cdots H_{n-1}^{m_{n-1}})=0
$$
establish a gradation in $U_q({\rm sl}_n)$ (see [7], subsection 6.1.5).

Let us introduce the elements
$$
{\tilde I}_{j,j-1}=F_{j-1}-qq^{-H_{j-1}}E_{j-1},\ \ \ \  j=2,3,\cdots ,n,
$$
of $U_q({\rm sl}_n)$. It is proved in [12] that there exists the
algebra homomorphism $\varphi : U'_q({\rm so}_n)\to U_q({\rm sl}_n)$
uniquely determined by the relations
$\varphi (I_{i+1,i})={\tilde I}_{i+1,i}$, $i=1,2,\cdots$, $n-1$.
The following theorem states that this homomorphism is an
isomorphism.
\medskip

\noindent
{\bf Theorem 3.} {\it The homomorphism
$\varphi : U'_q({\rm so}_n)\to U_q({\rm sl}_n)$
determined by the relations
$\varphi (I_{i+1,i})={\tilde I}_{i+1,i}$, $i=1,2,\cdots ,n-1$, is an
isomorphism of $U'_q({\rm so}_n)$ to $U_q({\rm sl}_n)$.}
\medskip

\noindent
{\sl Proof.}
In [22] the authors of that paper state that this homomorphism is
an isomorphism and say that it can be proved by means of the Diamond Lemma.
However, we could not restore their proof and found another one.
It use the above Poincar\'e--Birkhoff--Witt theorem for the
algebra $U'_q({\rm so}_n)$.
Namely, we use the explicit expressions from [22] for the
elements ${\tilde I}_{ij}\equiv \varphi (I_{ij})
\in U_q({\rm sl}_n)$ in terms of the elements of
the $L$-functionals of the quantum algebra $U_q({\rm sl}_n)$:
$$
{\tilde I}_{ji}=(q-q^{-1})^{-1}c_iK_{ji}^-,\ \ \  j>i, \eqno (17)
$$
where $c_i$ is equal to $q^s$ with an appropriate $s\in {\Bbb Z}$ and
$$
K^-\equiv (K_{ji}^-)_{i,j=1}^n =(L^+)^tJL^-. \eqno (18)
$$
Here $J={\rm diag}\, (q^{n-1},q^{n-2},\cdots ,1)$ and explicit expressions
for matrix elements $l^+_{ij}$ and $l^-_{ij}$ of the matrices $L^+$ and
$L^-$ are given by formulas from [12] (see also [7],
subsection 8.5.2). In particular,
$l^+_{ij}=l^-_{ji}=0$ if $i>j$ and $l^+_{ij}$ (resp.
$l^-_{ji}$) is expressed in terms of $E_{\beta _{ij}}$ (in terms of
$F_{\beta _{ij}}$) if $i<j$. We have
$$
{\rm deg}\, l^+_{ij}=\beta _{ij},\ \ \
{\rm deg}\, l^-_{ij}=-\beta _{ij}.
$$
The elements $l^\pm _{jj}$ belong to the subalgebra ${\frak H}$
generated by $K_1,K_2,\cdots ,K_{n-1}$. By (18) for $j>i$ we obtain
$$
K^-_{ji}=\sum _{s=i}^j c'_sl^+_{sj}l^-_{si}, \ \ \ \ j>i,
\eqno (19)
$$
where $c'_s=(q-q^{-1})q^r$ with an appropriate $r\in {\Bbb Z}$.
The summands in (19) have different degrees and the lowest degree
has the only summand $c'_jl^+_{jj}l^-_{ji}$.

Let $a$ be a basis element
${I_{21}^-}^{m_{21}}{I_{31}^-}^{m_{31}}\cdots
{I_{n,n-1}^-}^{m_{n,n-1}}$
of the algebra $U'_q({\rm so}_n)$
from Theorem 2. Then
$$
\varphi (a)=
({\tilde I}^-_{21})^{m_{21}}({\tilde I}^-_{31})^{m_{31}}\cdots
({\tilde I}^-_{n,n-1})^{m_{n,n-1}}.
$$
Substituting here expressions for ${\tilde I}^-_{ji}$ from
formulas (17) and (19), we
obtain $\varphi (a)$ in form of a sum with a single summand
of the lowest degree. This summand of lowest degree is
$c'F_{\beta _{21}}^{m_{21}}F_{\beta _{31}}^{m_{31}}\cdots
F_{\beta_{n,n-1}}^{m_{n,n-1}}$ with nonvanishing coefficient $c'$.
The expression at $c'$
is just the basis element of ${\frak N}^-\subset U_q({\rm sl}_n)$
corresponding under the mapping ${\cal T}$ to the basis element
$a$ of $U'_q({\rm so}_n)$.

Similarly, if an element $a\in U'_q({\rm so}_n)$
is a linear combination of the basis elements
${I_{21}^-}^{m_{21}}{I_{31}^-}^{m_{31}}\cdots
{I_{n,n-1}^-}^{m_{n,n-1}}$
from Theorem 2, then we substitute into $\varphi (a)$ expressions
(19) for each $K^-_{ji}$. As a result, we express $\varphi (a)$ in
form af a sum, containing the same linear combination of products
$c'F_{\beta_{21}}^{m_{21}}F_{\beta_{31}}^{m_{31}}\cdots
F_{\beta_{n,n-1}}^{m_{n,n-1}}$. This linear combination contains
a subsum of a (fixed) lowest degree and this subsum cannot be
concelled with other summands in $\varphi (a)$. Therefore, $\varphi (a)
\ne 0$ and
$\varphi$ is an isomorphism from $U'_q({\rm so}_n)$ to $U_q({\rm sl}_n)$.
Theorem is proved.
\medskip

This theorem has an important corollary.
\medskip

\noindent
{\bf Corollary.}
{\it Finite dimensional irreducible representations
of $U'_q({\rm so}_n)$ separate elements of this algebra, that is,
for any $a\in U'_q({\rm so}_n)$ there exists a finite dimensional
irreducible representation $T$ of $U'_q({\rm so}_n)$ such that
$T(a)\ne 0$.}
\medskip

\noindent
{\sl Proof.} If $q$ is not a root of unity, then the assertion of the
theorem follows from Theorem 3 and from the theorem on
separation of elements of the algebra $U_q({\rm sl}_n)$ by its
representations (see subsection 7.1.5 in [7]) if to take into account
the fact (proved in [22]) that a restriction of
any finite dimensional irreducible representation
of $U_q({\rm sl}_n)$ onto the subalgebra $U'_q({\rm so}_n)$
decomposes into a direct sum of its irreducible representations.

Let now $q$ be a root of unity, that is
$q^k=1$. We denote by $N$ a positive integer such that every
irreducible representation of the quantum algebra $U_q({\rm sl}_n)$
has dimension less than $N$. Let $a$ be any nonvanishing element
of $U'_q({\rm so}_n)$. Then there exists an irreducible representation
$T$ of $U_q({\rm sl}_n)$ such that $T(\varphi (a^N))\ne 0$. (Note that
$a^N\ne 0$ since $U_q({\rm sl}_n)$ has no divisors of zero.) Let
${\tilde T}$ be the restriction of $T$ to the subalgebra
$U'_q({\rm so}_n)$. Then ${\tilde T}$ is reducible. For simplicity we
suppose that ${\tilde T}$ contains only two irreducible representations
of $U'_q({\rm so}_n)$. (If ${\tilde T}$ contains more irreducible
constituents, then the proof is the same as for two ones.) Generally
speaking, ${\tilde T}$ is not completely reducible, that is in some
basis the representation ${\tilde T}$ is of the form
$$
\left( \matrix{ T_1 & *\cr
                 0  & T_2} \right) ,
$$
where $T_1$ and $T_2$ are irreducible representations of $U'_q({\rm so}_n)$.
Since ${\tilde T}(a^N)\ne 0$, then ${\tilde T}(a)\ne 0$ and
$$
{\tilde T}(a)=
\left( \matrix{ T_1(a) & *\cr
                 0  & T_2(a)} \right) .
$$
If $T_1(a)\ne 0$ or $T_2(a)\ne 0$, then irreducible representations of
$U'_q({\rm so}_n)$ separate the element $a$. Let
$T_1(a)=0$ and $T_2(a)=0$. Then
$$
{\tilde T}(a)=
\left( \matrix{ 0 & *\cr
                 0  & 0} \right) .
$$
In this case ${\tilde T}(a)$ is a nilpotent matrix and
${\tilde T}(a)^N={\tilde T}(a^N)=0$. This contradict the assupmtion that
${\tilde T}(a^N)\ne 0$. Therefore, the case $T_1(a)=0$ and $T_2(a)=0$
is not possible. Corollary is proved.
\bigskip

\centerline{\sc 4. Finite dimensionality of representations}
\medskip

\noindent
Everywhere below we assume, if other is not stated, that $q$ is
a root of unity. Moreover, we consider that $q^k=1$ and $k$ is an
odd integer.

Below we shall need an information on the center of the algebra
$U'_q({\rm so}_n)$.
Central elements of the algebra $U'_q({\rm so}_n)$ for
any value of $q$ are found in [17] and [22]. They are given in the form of
homogeneous polynomials of elements of $U'_q({\rm so}_n)$.
If $q$ is a root of unity, then (as in the case of Drinfeld--Jimbo
quantum algebras)
there are additional central elements of $U'_q({\rm so}_n)$ which are
given by the following theorem, proved in [18].
\medskip

\noindent
{\bf Theorem 4.} {\it Let $q^k=1$ for $k \in {\Bbb N}$ and $q^j \not=1$
for $0<j< k$. Then the elements
$$
C^{(k)}(I^+_{rl})=
\sum_{j=0}^{\{ (k-1)/2\}} \left( \matrix{k-j\cr j}\right) \frac{1}{k-j}
\Bigl( \frac{i}{q-q^{-1}} \Bigr)^{2j} {I^+_{rl}}^{k-2j},\ \ \ r>l,
\eqno(20)
$$
where ${\{ (k-1)/2\}}$ is the integral part of the number $(k-1)/2$,
belong to the center of $U'_q({\rm so}_n)$.}
\medskip

It is well-known that a Drinfeld--Jimbo algebra $U_q(g)$ for $q$ a root
of unity ($q^k=1$) is a finite dimensional vector space
over the center of $U_q(g)$. The same assertion is true for the
algebra $U'_q({\rm so}_n)$. In fact, by Theorem 4 any element
$(I^+_{ij})^s$, $s\ge k$, can be reduced to a linear combination of
$(I^+_{ij})^r$, $r<k$, with coefficients from the center ${\cal C}$
of $U'_q({\rm so}_n)$. Now our assertion follows from this sentence
and from Poincar\'e--Birkhoff--Witt theorem for $U'_q({\rm so}_n)$.
\medskip

\noindent
{\bf Theorem 5.} {\it If $q$ is a root of unity, then any irreducible
representation of $U'_q({\rm so}_n)$ is finite dimensional.}
\medskip

\noindent
{\sl Proof.} Let $q$ be a root of unity, that is $q^k=1$.
Let $T$ be an irreducible representation of $U'_q({\rm so}_n)$
on a vector space ${\cal V}$.
Then $T$ maps central elements into scalar operators. Since the linear
space $U'_q({\rm so}_n)$ is finite dimensional over the center ${\cal C}$
with the basis
${I_{21}^+}^{m_{21}}{I_{31}^+}^{m_{31}}\cdots
{I_{n,n-1}^+}^{m_{n,n-1}}$, $m_{ij}<k$,
then
for any $a\in U'_q({\rm so}_n)$ we have $T(a)=\sum _{ m_{ij}<k}
c_{\{ m_{ij}\} }T({I_{21}^+}^{m_{21}}{I_{31}^+}^{m_{31}}\cdots
{I_{n,n-1}^+}^{m_{n,n-1}})$, where
$c_{\{ m_{ij}\} }$ are numerical coefficients.
Hence, if ${\bf v}$ is a nonzero vector of the
representation space ${\cal V}$, then $T(U'_q({\rm so}_n)){\bf v}={\cal V}$
since $T$ is an irreducible representation. Since $T(a)$ is of the above
form for any $a\in U'_q({\rm so}_n)$, then
${\cal V}$ is finite dimensional. Theorem is proved.
\medskip

It follows from this proof that there exists a fixed positive integer
$r$ such that dimension of any irreducible representation of
$U'_q({\rm so}_n)$ at $q$ a root of unity does not exceed $r$. Of course,
the number $r$ depends on $k$ (recall that $k$ is defined by $q^k=1$).
\bigskip

\centerline{\sc 5. Cyclic representations
at $q$ a root of unity}
\medskip

\noindent
Taking into account Theorem 5, below under studying irreducible
representations of $U'_q({\rm so}_n)$ at $q$ a root of unity
we consider only its finite dimensional representations.

If $q$ is not a root of unity, there
exists two types of such representations:
\medskip

(a) representations of
the classical type (at $q\to 1$ they give the corresponding
finite dimensional irreducible representations of
the Lie algebra ${\rm so}_n$);

(b) representations of the nonclassical
type (they do not admit the limit $q\to 1$ since in this point the
representation operators are singular). These representations are
described in [14].
\medskip

Let us consider irreducible representations of $U'_q({\rm so}_n)$ for
$q$ a root of unity
($q^k=1$ and $k$ is a smallest positive integer with this property).
We also assume that $k$ is odd. If $k$ would be even, then almost
all below reasoning is true, if to replace $k$ by $k'=k/2$
(as in the case of irreducible representations of the quantum algebra
$U_q({\rm sl}_2)$ for $q$ a root of unity in [7], chapter 3).

We fix complex numbers
$m_{1,n}, m_{2,n},..., m_{\left \{ {n/2}\right \} ,n}$
(here $\{ {n/2}\}$ denotes an integral part
of ${n/2}$) and $c_{ij}$, $h_{ij}$, $j=2,3,\cdots ,n-1$,
$i=1,2,\cdots , \{ {j/2}\}$ such that no of the numbers
$$
m_{in},\ \ h_{ij},\ \
h_{ij}-h_{sj},\ \ h_{ij}-h_{s,j\pm 1},\ \  h_{ij}+h_{sj},\ \
h_{ij}+h_{s,j\pm 1},\ \ h_{b,n-1}-m_{sn},\ \ h_{b,n-1}+m_{sn}
$$
belongs to $\frac 12 {\Bbb Z}$. (We also suppose that $c_{ij}\ne 0$.)
The set of these numbers will be denoted by $\omega$:
$$
\omega =\{ {\bf m}_n, {\bf c}_{n-1}, {\bf h}_{n-1},\cdots ,
{\bf c}_2, {\bf h}_2\} ,
$$
where ${\bf m}_n$ is the set of the numbers
$m_{1,n}, m_{2,n},..., m_{\left \{ {n/2}\right \} ,n}$, and
${\bf c}_{j}$ and ${\bf h}_{j}$ are the sets of numbers
$c_{ij}$, $i=1,2,\cdots ,\{ j/2\}$, and
$h_{ij}$, $i=1,2,\cdots ,\{ j/2\}$, respectively.
(Thus, $\omega$ contains $r={\rm dim}\, {\rm so}_n$ complex numbers.)
Let $V$ be a complex vector space with a basis labelled by the
tableaux
$$
  \{\xi_{n} \}
\equiv \left\{ \matrix{ {\bf m}_{n} \cr {\bf m}_{n-1} \cr \dots
\cr {\bf m}_{2}  }
 \right\}
\equiv \{ {\bf m}_{n},\xi_{n-1}\}\equiv \{{\bf m}_{n} ,
{\bf m}_{n-1} ,\xi_{n-2}\} ,
                                                       \eqno(21)
$$
where the set of numbers ${\bf m}_{n}$
consists of $\{ {n/2}\}$ numbers $m_{1,n}, m_{2,n},\cdots , m_{\left
\{ {n/2}\right \} ,n}$ given above, and for each $s=2,3,\cdots ,n-1$,
${\bf m}_{s}$ is a set of numbers $m_{1,s}, \cdots , m_{\left
\{ {s/2}\right \} ,s}$ and each $m_{i,s}$ runs independently the
values $h_{i,s}, h_{i,s}+1,\cdots , h_{i,s}+k-1$.
Thus, ${\rm dim}\, V$ coincides with $k^N$, where $N$ is the number of
positive roots of ${\rm so}_n$.
It is convenient to use for the numbers $m_{i,s}$, $s=2,3,\cdots ,n$,
the so-called $l$-coordinates
$$
l_{j,2p+1}=m_{j,2p+1}+p-j+1,  \qquad
                      l_{j,2p}=m_{j,2p}+p-j .        \eqno (22)
$$

To the set of numbers $\omega$ there corresponds the
irreducible finite dimensional representation $T_\omega$ of the algebra
$U'_q({\rm so}_n)$.
The operators $T_{\omega}(I_{2p+1,2p})$ of the representation
$T_{\omega}$ act upon the basis elements, labelled by (21), by the formula
$$
T_{\omega}(I_{2p+1,2p})
| \xi_n\rangle =
\sum^p_{j=1} c_{j,2p} \frac{ A^j_{2p}(\xi_n)}
{q^{l_{j,2p}}+q^{-l_{j,2p}} } \,
            \vert (\xi_n)^{+j}_{2p}\rangle -
\sum^p_{j=1} c^{-1}_{j,2p} \frac{A^j_{2p}((\xi_n)^{-j}_{2p})}
{q^{l_{j,2p}}+q^{-l_{j,2p}}} \,
|(\xi_n)^{-j}_{2p}\rangle                         \eqno(23)
$$
and the operators $T_{\omega}(I_{2p,2p-1})$ of the representation
$T_{\omega}$ act as
$$
T_{\omega}(I_{2p,2p-1})\vert \xi_n\rangle=
\sum^{p-1}_{j=1}c_{j,2p-1}  \frac{B^j_{2p-1}(\xi_n)}
{[2 l_{j,2p-1}-1][l_{j,2p-1}]} \,
\vert (\xi_n)^{+j}_{2p-1} \rangle -
$$
$$
-\sum^{p-1}_{j=1}c^{-1}_{j,j,2p-1} \frac {B^j_{2p-1}((\xi_n)^{-j}_{2p-1})}
{[2 l_{j,2p-1}-1][l_{j,2p-1}-1]} \,
\vert (\xi_n)^{-j}_{2p-1}\rangle
+ {\rm i}\, C_{2p-1}(\xi_n) \,
\vert \xi_n \rangle ,                                 \eqno(24)
$$
where numbers in square brackets mean $q$-numbers:
$$
[b]:= \frac{q^b-q^{-b}}{q-q^{-1}}.
$$
In these formulas, $(\xi_n)^{\pm j}_{s}$ means the tableau (21)
in which $j$-th component $m_{j,s}$ in ${\bf m}_s$ is replaced
by $m_{j,s}\pm 1$. If $m_{j,s}+1=h_{j,s}+k$ (resp.
$m_{j,s}-1=h_{j,s}-1$), then we set $m_{j,s}+1=h_{j,s}$ (resp.
$m_{j,s}-1=h_{j,s}+k-1$). The coefficients
$A^j_{2p},  $ $B^j_{2p-1},$ $C_{2p-1}$ in (23) and (24) are given
by the expressions
$$ \leqno
A^j_{2p}(\xi_n) =
$$  $$
=\left( \frac{\prod_{i=1}^p [l_{i,2p+1}+l_{j,2p}]
[l_{i,2p+1}{-}l_{j,2p}{-}1]
\prod_{i=1}^{p-1} [l_{i,2p-1}+l_{j,2p}] [l_{i,2p-1}{-}l_{j,2p}{-}1]}
{\prod_{i\ne j}^p [l_{i,2p}+l_{j,2p}][l_{i,2p}-l_{j,2p}]
[l_{i,2p}+l_{j,2p}+1][l_{i,2p}-l_{j,2p}-1]} \right)^{1/2} ,
$$
$$\leqno
B^j_{2p-1}(\xi_n)=
$$   $$
=\left( \frac{\prod_{i=1}^p
[l_{i,2p}+l_{j,2p-1}] [l_{i,2p}-l_{j,2p-1}] \prod_{i=1}^{p-1}
[l_{i,2p-2}+l_{j,2p-1}] [l_{i,2p-2}-l_{j,2p-1}]}
{\prod_{i\ne j}^{p-1}
[l_{i,2p-1}{+}l_{j,2p-1}][l_{i,2p-1}{-}l_{j,2p{-}1}]
[l_{i,2p-1}{+}l_{j,2p-1}{-}1][l_{i,2p-1}{-}l_{j,2p-1}{-}1]}
\right) ^{1/2} ,
$$
$$
C_{2p-1}(\xi_n) =\frac{ \prod_{s=1}^p [ l_{s,2p} ]
\prod_{s=1}^{p-1} [ l_{s,2p-2} ]}
{\prod_{s=1}^{p-1} [l_{s,2p-1}] [l_{s,2p-1} - 1] } .
$$

The fact that the operators $T_\omega (I_{j,j-1})$, given above,
satisfy the defining relations (1)--(3) is proved in the same
way as in the case of irreducible representations of
$U'_q({\rm so}_n)$ when $q$ is not a root of unity in [13] and we omit
these rather long calculations.

As in the case of finite dimensional irreducible representations of the
Lie algebra ${\rm so}_n$, the form of the basis elements of the above
representation space $V$ and the formulas for the operators $T_\omega
(I_{j,j-1})$ allow us to decompose the restriction of the representation
$T_\omega$,
$\omega =\{ {\bf m}_n, {\bf c}_{n-1}, {\bf h}_{n-1},\cdots ,
{\bf c}_2, {\bf h}_2\}$, to the subalgebra $U'_q({\rm so}_{n-1})$.
We have
$$
T_\omega \bigl| _{U'_q({\rm so}_{n-1})}=\bigoplus \nolimits_{\omega _{n-1}}
T_{\omega _{n-1}} , \eqno (25)
$$
where
$\omega_{n-1} =\{ {\bf m}_{n-1}, {\bf c}_{n-2}, {\bf h}_{n-2},\cdots ,
{\bf c}_2, {\bf h}_2\}$ and ${\bf m}_{n-1}$ runs over the vectors
$$
(h_{1,n-1}+a_1, h_{2,n-1}+a_2,\cdots ,h_{s,n-1}+a_s),\ \ \
s=\{ (n-1)/2\} ,\ \  a_j=0,1,2,\cdots ,k-1,
$$
and ${\bf c}_j$ and ${\bf h}_j$ are such as in $\omega$.
\medskip

\noindent
{\bf Theorem 6.} {\it Representations $T_\omega$ with the domain of values
of representation parameters, as described above, are irreducible.}
\medskip

\noindent
{\sl Proof.}
The proof will be carried out by induction. We shall show even more:
For every algebra $U'_q({\rm so}_r)$, $r=2,3,\cdots$,
the representations $T_{\omega}$ are irreducible, and two
representations $T_{\omega}$ and $T_{\omega'}$ ($\omega\ne\omega'$)
with $\omega=\{{\bf m}_r, {\bf c}_{r-1}, {\bf h}_{r-1}, \cdots,
{\bf c}_{2}, {\bf h}_{2}\}$ and
$\omega'=\{{\bf m}'_r, {\bf c}_{r-1}, {\bf h}_{r-1}, \cdots,
{\bf c}_{2}, {\bf h}_{2}\}$, such that the numbers $l_{i,r}-l'_{i,r}$,
$1\le i\le \{r/2\}$, are integers and $-k<l_{i,r}-l'_{i,r}<k$,
are nonequivalent.

For the algebra $U'_q({\rm so}_2)$ these statements are true. In fact,
the irreducible representations of $U'_q({\rm so}_2)$ are
one-dimensional and nonequivalence condition means that
${\rm i}[m_{12}]\ne {\rm i}[m'_{12}]$ if $m_{12}-m'_{12}$ is an integer
such that $-k<m_{12}-m'_{12}<k$.
This is true due to the properties of $q$-numbers $[m]$ for $q^k=1$.

Now we assume that the above statements on irreducibility and
nonequivalence are true for the representations $T_\omega$ of the algebra
$U'_q({\rm so}_{n-1})$ and prove that they are true for the representations
of $U'_q({\rm so}_n)$.

First we show that the representations $T_{\omega}$
of $U'_q({\rm so}_n)$ are irreducible.
For this end, we use the decomposition (25):
$$
T_\omega\bigr| _ {U'_q({\rm so}_{n-1})} = T_1\oplus T_2\oplus\cdots\ ,
$$
where $T_i$ are the corresponding representations $T_{\omega _{n-1}}$.
According to the induction assumption, all the representations $T_i$
in this decomposition are irreducible and pairwise nonequivalent.
The corresponding representation space $V_\omega$ decomposes as
$$
V_\omega=V_1\oplus V_2\oplus\cdots\ , \eqno (26)
$$
where $V_i$ are representation
spaces for $T_i$. Let $T_\omega$ be reducible. This means that there
exists a proper subspace $V$ in $V_\omega$ which is
invariant under $U'_q({\rm so}_n)$. Let us show that
$V$ is a direct sum of some of the subspaces $V_i$. In fact, if
${\bf x}\in V$, then
${\bf x}={\bf x}_1+{\bf x}_2+\ldots+{\bf x}_r$, where
${\bf x}_1 \in V_{i_1}$, ${\bf x}_2 \in V_{i_2}$, $\cdots$,
${\bf x}_r \in V_{i_r}$. First suppose that $r=2$ and ${\bf x}=
{\bf x}_1+{\bf x}_2$, ${\bf x}_1 \ne 0$, ${\bf x}_2 \ne0$.
Then $T_\omega (U'_q({\rm so}_{n-1})){\bf x}$ is a subspace
$V'$ of $V$, invariant with respect to
$U'_q({\rm so}_{n-1})$. Since $V_{i_1}$ and $V_{i_2}$ are irreducible
for $U'_q({\rm so}_{n-1})$ and the representations $T_{i_1}$ and $T_{i_2}$
are nonequivalent, then $V'$ coincides with one of the spaces
$V_{i_1}$, $V_{i_2}$, $V_{i_1}\oplus V_{i_2}$. In any case,
$V_{i_1}\subset V'$, $V_{i_2}\subset V'$ and
$V'=V_{i_1}\oplus V_{i_2}$. This means that
$V_{i_1}\oplus V_{i_2}\subset V$. It is proved similarly that if $r>2$,
then $V_{i_1}\oplus \cdots \oplus V_{i_r}\subset V$.
Thus, if some vector ${\bf x}\in V$ has nonzero component in some
subspace $V_i$ from the decomposition (26) for
$V_\omega$, then $V_i\subset V$.

Now it follows from the above formulas for the operators
$T_\omega(I_{n,n-1})$ that the action of the operator
$T_\omega(I^s_{n,n-1})$, $s=0,1,2,\cdots$, upon
arbitrary vector ${\bf x}\in V_i$ gives vectors containing nonzero
components of any subspace $V_j$ from (26). Therefore,
any subspace $V_j$ from (26) belongs to $V$. This means that
$V=V_\omega$ and the representation $T_\omega$ is irreducible.

Let us show nonequivalence of two
representations $T_{\omega}$ and $T_{\omega'}$, $\omega\ne\omega'$,
with
$\omega=\{{\bf m}_n, {\bf c}_{n-1}, {\bf h}_{n-1}, \cdots,
{\bf c}_{2}, {\bf h}_{2}\}$ and
$\omega'=\{{\bf m}'_n, {\bf c}_{n-1}, {\bf h}_{n-1}, \cdots,
{\bf c}_{2}, {\bf h}_{2}\}$, where all the numbers $l_{i,n}-l'_{i,n}$,
$1\le i\le \{n/2\}$, satisfy the conditions formulated in the beginning
of this proof. Assume that these representations are equivalent and
show that this leads to contradiction.
By the definition, two irreducible representations
$T_\omega$ and $T_{\omega'}$ are equivalent if there exists
nondegenerate operator $A$ such that
$$
A T_\omega(a) A^{-1} = T_{\omega'}(a),\ \ \ \ a\in U'_q({\rm so}_n).
\eqno(27)
$$
The representations $T_\omega$ and $T_{\omega'}$ have the same
decomposition into irreducible (and pairwise nonequivalent)
representations under restriction
onto $U'_q({\rm so}_{n-1})$. (Note that
$T_\omega$ and $T_{\omega'}$ are defined on the same space.)
Due to Schur lemma, this means that the operator $A$ is a direct sum of
operators $A_i$ acting on the subspaces $V_i$ from (26)
and each $A_i$ is multiple to the unit operator.

Let us consider the case $n=2p$. Putting $a=I_{2p,2p-1}$ in (27) and
writing this relation in matrix form in the basis (21) we
obtain the equalities
$$
B^j_{2p-1}(\xi_n) \frac{a_{{\bf m}_{n-1}^{+j}}}
{a_{{\bf m}_{n-1}}}= B^j_{2p-1}(\xi'_n) , \eqno(28)
$$
$$
B^j_{2p-1}(\xi_n^{-j}) \frac{a_{{\bf m}_{n-1}^{-j}}}
{a_{{\bf m}_{n-1}}}= B^j_{2p-1}(\xi^{'-j}_n) ,
$$
where $\xi_n$ differs from $\xi '_n$ by the replacement
$m_{j,2p}\to m'_{j,2p}$ and $a_{{\bf m}_{n-1}}$ are matrix elements of
the matrix $A$.
Replacing $\xi_n$ by $\xi_n^{+j}$ in the last relation we have
$$
B^j_{2p-1}(\xi_n) \frac{a_{{\bf m}_{n-1}}}
{a_{{\bf m}_{n-1}^{+j}}}= B^j_{2p-1}(\xi'_n) . \eqno(29)
$$
Eliminating all $a_{{\bf m}_{n-1}}$ from (28) and (29)
we obtain
$$
(B^j_{2p-1}(\xi_n))^2=(B^j_{2p-1}(\xi'_n))^2.    \eqno(30)
$$
Using the explicit form of the coefficients $B^j_{2p-1}$
and the identity $[x-y][x+y]=[x]^2-[y]^2$, we derive that
$$
\prod_{i=1}^p ([l_{i,2p}]^2-[l_{j,2p-1}]^2)=
\prod_{i=1}^p ([l'_{i,2p}]^2-[l_{j,2p-1}]^2). \eqno(31)
$$

Let us consider two polynomials
$\prod_{i=1}^p ([l_{i,2p}]^2-x)$ and
$\prod_{i=1}^p ([l'_{i,2p}]^2-x)$ of variable $x$.
They coincide identically if
they coincide on $p$ different values of $x$. Equalities (31)
give such values. Indeed, for each fixed $j$ ($1\le j\le p-1$),
$[l_{j,2p-1}]^2$ takes $k$ different values, and sets of values
obtained for different $j$ do not intersect. Thus, we have
coincidence of the above polynomials on $k (p-1)$ different
points. Therefore, they coincide identically, and the sets of their
zeros $\{[l_{i,2p}]^2\}$ and $\{[l'_{i,2p}]^2\}$, $1\le i\le p$,
also must coincide. But the latter is impossible due to the
properties of $q$-numbers and due to
conditions for $\{l_{i,2p}\}$ and $\{l'_{i,2p}\}$ (recall that
$l_{i,2p},l'_{i,2p}\not\in \frac 12 {\Bbb Z}$).

A nonequivalence of the representations $T_\omega$ in the case of the
algebra $U'_q({\rm so}_{2p+1})$ is proved analogously.
We only note that in this case we obtain the relations
$(A^j_{2p}(\xi_n))^2=(A^j_{2p}(\xi'_n))^2$
(instead of relations (30)) reducing (due to
the identity $[x][y]=[(x+y)/2]^2-[(x-y)/2]^2$) to
$$
\prod_{i=1}^p ([l_{i,2p+1}-1/2]^2-[l_{j,2p}+1/2]^2)=
\prod_{i=1}^p ([l'_{i,2p+1}-1/2]^2-[l_{j,2p}+1/2]^2).
$$
Theorem is proved.
\medskip

There are equivalence relations in the set of irreducible
representations $T_\omega$. In order to extract a
subset of pairwise nonequivalent representations from the entire set,
we introduce some domains on the complex plane. The set
$$
D=\{x\in {\Bbb C}\; |\; |{\rm Re}\, x|<k/4 \ {\rm or}\
{\rm Re}\, x = -k/4,\ {\rm Im}\, x \le 0 \ {\rm or}\
{\rm Re}\, x = k/4,\ {\rm Im}\, x \ge 0\}
$$
is a maximal subset of ${\Bbb C}$
such that for all $x,y\in D$, $x \ne y$, we have $[x]\ne[y]$.
The set
$$
D^\pm=\{x\in {\Bbb C}\; |\; 0<{\rm Re}\, x<k/4 \ {\rm or}\
{\rm Re}\, x = 0,\ {\rm Im}\, x \ge 0 \ {\rm or}\
{\rm Re}\, x = k/4,\ {\rm Im}\, x \ge 0\}
$$
is a maximal subset of ${\Bbb C}$
such that for all $x,y\in D^\pm$, $x \ne y$, we have $[x]\ne\pm[y]$.
We need also the sets
$$
D_h=\{x\in {\Bbb C}\; |\; |{\rm Re}\, x|<1/4 \ {\rm or}\
{\rm Re}\, x = -1/4,\ {\rm Im}\, x \le 0 \ {\rm or}\
{\rm Re}\, x = 1/4,\ {\rm Im}\, x \ge 0\},
$$
$$
D^\pm_h=\{x\in {\Bbb C}\; |\; 0<{\rm Re}\, x<1/4 \ {\rm or}\
{\rm Re}\, x = 0,\ {\rm Im}\, x \ge 0 \ {\rm or}\
{\rm Re}\, x = 1/4,\ {\rm Im}\, x \ge 0\}.
$$

We introduce an ordering in the set $D^\pm$ (resp. $D^\pm_h$)
as follows:
we say that $x\succ y$, $x,y\in D^\pm$ (resp. $x,y\in D^\pm_h$) if
either ${\rm Re}\, x > {\rm Re}\, y$ or both ${\rm Re}\, x = {\rm Re}\, y$
and ${\rm Im}\, x > {\rm Im}\, y$.

We say that the set of complex numbers
${\bf l}_{2p}=(l_{1,2p},l_{2,2p},\cdots,l_{p,2p})$ is {\it
dominant} if $l_{1,2p}$, $l_{2,2p}, \cdots$, $l_{p-1,2p}\in D^\pm$,
$l_{p,2p}\in D$, and $l_{1,2p}\succ l_{2,2p} \succ\cdots\succ
l_{p-1,2p}\succ l^*_{p,2p}$, where
$l^*_{p,2p}=l_{p,2p}$ if $l_{p,2p}\in D^\pm$ and
$l^*_{p,2p}=-l_{p,2p}$ if $l_{p,2p}\not\in D^\pm$.

The notion of {\it dominance} for
the set ${\bf h}_{2p}=(h_{1,2p},h_{2,2p},\cdots,
h_{p,2p})\in {\Bbb C}^p$ is introduced by the replacements
$l_{i,2p}\to h_{i,2p}$, $D\to D_h$ and $D^\pm\to D^\pm_h$ in the
previous definition.

We say that the set of complex numbers
${\bf l}_{2p+1}=(l_{1,2p+1},l_{2,2p+1},\cdots,l_{p,2p+1})$ is {\it
dominant} if $l_{1,2p+1}$, $l_{2,2p+1},\cdots,l_{p,2p+1}\in D^\pm$
and $l_{1,2p+1}\succ l_{2,2p+1} \succ\cdots\succ
l_{p,2p+1}$.

The notion of {\it dominance} for
the set of complex numbers ${\bf h}_{2p+1}=(h_{1,2p+1},h_{2,2p+1}$,
$\cdots,h_{p,2p+1})$ is introduced by the replacements
$l_{i,2p+1}\to h_{i,2p+1}$ and $D^\pm\to D^\pm_h$ in the
previous definition.

We say that
$\omega=\{{\bf m}_n, {\bf c}_{n-1}, {\bf h}_{n-1}, \cdots, {\bf c}_2, {\bf
h}_2\}$ is {\it dominant} if every of the sets ${\bf l}_n$,
${\bf h}_{n-1}, \cdots$, ${\bf h}_2$ is dominant
and if $0\le \mathop{\rm Arg} c_{ij} < 2\pi/k,$  $j=2,3,\cdots ,n-1$;
$i=1,2,\cdots ,\{ j/2\}$.
\medskip

\noindent
{\bf Theorem 7.} {\it The representations $T_\omega$ of
$U'_q({\rm so}_n)$ with dominant
$\omega$ are pairwise nonequivalent. Any irreducible representation
$T_{\omega '}$ is equivalent to some representation $T_\omega$
with dominant $\omega$.}
\medskip

This theorem is proved by using the relation similar to relation (27)
and the decomposition (25) for the restriction of the representations
$T_\omega$ to the subalgebra $U'_q({\rm so}_{n-1})$.
\bigskip

\centerline{\sc 6. Partially cyclic representations
at $q$ a root of unity}
\medskip

\noindent
The representations of the previous section constitute a main class of
irreducible representations of $U'_q({\rm so}_n)$ at $q$ a root of
unity. There are many other classes of irreducible representations
which are given by smaller number of parameters and act on linear spaces
of smaller dimension. They are obtained from the above
representations $T_\omega$ if the representation parameters take the
values excluded in the previous section. For these values of
parameters the corresponding representations become reducible and
their irreducible constituents constitute new classes of
irreducible representations of $U'_q({\rm so}_n)$. Let us give some
of these classes of representations.

Let us fix an integer $i$ such that $1\le i<\{ n/2\} $, where
$\{ n/2\} $ is an integral part of $n/2$.
We fix complex numbers
$m_{1n}, m_{2n},\cdots$, $m_{in}$ and $c_{sj}$, $h_{sj}$,
$j=2,3,\cdots ,n-1$, $s=1,2,\cdots , {\rm max}\, \{ i,\{ j/2\} \}$, such
that
$c_{sj}\ne 0$ and no of the numbers
$$
m_{rn},\ \ h_{rj},\ \
h_{rj}-h_{sj},\ \ h_{rj}-h_{s,j\pm 1},\ \  h_{rj}+h_{sj},\ \
h_{rj}+h_{s,j\pm 1},\ \ h_{r,n-1}-m_{sn},\ \ h_{r,n-1}+m_{sn}
$$
belongs to $\frac 12 {\Bbb Z}$, and the
numbers $m_{i+1,n},m_{i+2,n},\cdots ,m_{\{ n/2\} ,n}$, which are all
integral or all half-integral and such that
$$
m_{i+1,n}\ge m_{i+2,n}\ge \cdots \ge |m_{ln}|\ \ \
{\rm if}\ \ \ n=2l ,
$$   $$
m_{i+1,n}\ge m_{i+2,n}\ge \cdots \ge m_{ln}\ge 0 \ \ \
{\rm if}\ \ \ n=2l+1 ,
$$   $$
l_{i+1,n}+l_{\{ n/2\} ,n}<k,\ \ \ \ l_{i+1,n}-l_{\{ n/2\} ,n}<k,
$$
where $l_{jn}$ are determined by formula (22).
Let $\omega$ be the set of all these numbers $m_{jn}, c_{sj}, h_{sj}$.
(Note that the quantity of these numbers is less than in the set
$\omega$ in the previous section.)

Let $V$ be a complex vector space with the basis labelled by the tableaux
(21),
where the set of numbers ${\bf m}_{n}$
consists of $\{ {n/2}\}$ numbers $m_{1,n}, m_{2,n},\cdots , m_{\left
\{ {n/2}\right \} ,n}$ given above, and for each $s=2,3,\cdots ,n-1$,
${\bf m}_{s}$ is the set of numbers $m_{1,s}, \cdots , m_{\left
\{ {s/2}\right \} ,s}$ and each $m_{rs}$, $r\le i$, runs independently the
values $h_{r,s}, h_{r,s}+1,\cdots , h_{r,s}+k-1$ and numbers
$m_{rs}$, $r> i$, run all integers (if
$m_{i+1,n},m_{i+2,n},\cdots ,m_{\{ n/2\} ,n}$ are integers) or
half-integers (if $m_{i+1,n},m_{i+2,n},\cdots ,m_{\{ n/2\} ,n}$ are
half-integers) satisfying the same betweenness conditions as in
the case of irreducible representations of the classical type
for $q$ not a root of unity, that is, the conditions
$$
m_{i+1,2r+1}\ge m_{i+1,2r}\ge m_{i+2,2r+1}\ge m_{i+2,2r}\ge
\cdots \ge m_{r,2r+1}\ge m_{r,2r}\ge -m_{r,2r+1} ,
$$  $$
m_{i+1,2r}\ge m_{i+1,2r-1}\ge m_{i+2,2r}\ge
m_{i+2,2r-1}\ge \cdots \ge m_{r-1,2r-1}\ge |m_{r,2r}| .
$$

To the set of numbers $\omega$ there corresponds the irreducible
finite dimensional representation $T'_\omega$ of the algebra
$U'_q({\rm so}_n)$ with the operators $T'_\omega (I_{r+1,r})$
given by formulas (23) and (24) (with the same expressions for the
coefficients $A^j_{2p}$, $B^j_{2p-1}$ and $C_{2p-1}$) in which
$c_{j,s}=1$ for $j>i$.

The fact that the operators $T'_\omega(I_{j,j-1})$ satisfy the defining
relations (1)--(3) of the algebra $U'_q({\rm so}_n)$ and irreducibility
of the representations $T'_\omega$ are proved in the same way as in the
previous section.

A particular case of the above representations $T'_\omega$ constitute
the so called irreducible representations of minimal dimension. They
correspond to the case $i=1$ and $m_{2n}=m_{3n}=\cdots =m_{\{ n/2\} ,n}=0$.
These representations are given only by complex numbers
$m_{1,n}$, $h_{1,n-1},h_{1,n-2},\cdots ,h_{1,2}$,
$c_{n-1},c_{n-2},\cdots ,c_{2}$. We denote them by $T^{\rm min}_\omega$.
The space of such a representation
has a basis labelled by the tableaux
$$
\xi _n =
\left(
\matrix{ m_{1n} & & 0 & \cdots & \cdots &\cdots & 0\cr
          & m_{1,n-1} & 0 & \cdots &\cdots  & 0\cr
           & \cdots & \cdots & \cdots &\cdots &  \cr
           &  & m_{1,4} & & 0 &  \cr
           &  & &  m_{1,3}  &   &  \cr
           &  & & m_{1,2} &  &  &  } \right) \equiv |m_{1,n},m_{1,n-1},
\cdots ,m_{1,2}\rangle ,
$$
where each $m_{1r}$, $r<n$, runs independently the values
$h_{1r},h_{1r}+1,\cdots , h_{1r}+k-1$. The operators
$T^{\rm min}_\omega (I_{r+1,r})$ of these representations
act upon these basis vectors by the formula
$$
T^{\rm min}_\omega (I_{r+1,r})|\xi_n \rangle =
$$   $$
=c_r \left(
\frac{ [m_{r+1}+m_{r}+r-2][m_{r+1}-m_{r}]
 [m_{r}+m_{r-1}+r-3]}
{[2m_{r}+r-1][2m_{r}+r-3][m_{r}-m_{r-1}+1]^{-1}}
\right) ^{1/2} | (\xi _n)^{+1}_r\rangle -
$$   $$
-c_r^{-1}\left(
\frac{ [m_{r+1}+m_{r}+r-3][m_{r+1}-m_{r}+1]
 [m_{r}+m_{r-1}+r-4]}
{[2m_{r}+r-3][2m_{r}+r-5][m_{r}-m_{r-1}]^{-1}}
\right) ^{1/2} | (\xi _n)^{-1}_r\rangle
$$
(for convenience we replaced here $m_{1,j}$ by $m_j$).
These representations are given by $2n-3$ parameters and act on
$k^{n-1}$-dimensional vector spaces.
\bigskip

\centerline{\sc 7. Other irreducible representations
at $q$ a root of unity}
\medskip

\noindent
If $q$ is a root of unity, then there also exists a class of irreducible
representations of $U'_q({\rm so}_n)$ similar to the representations
of the nonclassical type of $U'_q({\rm so}_n)$ when $q$ is not a root
of unity. These representations are described as follows.

Let $i$ be a fixed integer such that $1\le i<\{ n/2\}$.
We fix complex numbers
$m_{1n}, m_{2n},\cdots$, $m_{in}$ and $c_{sj}$, $h_{sj}$,
$j=2,3,\cdots ,n-1$, $s=1,2,\cdots , {\rm max}\, \{ i,\{ j/2\} \}$,
satisfying the same conditions as in section 6, and the
numbers $m_{i+1,n},m_{i+2,n},\cdots ,m_{\{ n/2\} n}$, which are all
half-integral and such that
$$
m_{i+1,n}\ge m_{i+2,n}\ge \cdots \ge m_{\{ n/2\} n}\ge 1/2 ,
$$ $$
l_{i+1,n}+l_{\{ n/2\} ,n}<k,\ \ \ \ l_{i+1,n}-l_{\{ n/2\} ,n}<k.
$$
We also fix the set $\epsilon =(\epsilon _{2i+2},\epsilon _{2i+3},\cdots ,
\epsilon _n)$, $\epsilon _j=\pm 1$.
Let $\omega$ be the set of all these numbers $m_{jn}, c_{sj}, h_{sj}$ and
$\epsilon$.

Let $V$ be a complex vector space with a basis labelled by the tableaux
(21),
where the set of numbers ${\bf m}_{n}$
consists of $\{ {n/2}\}$ numbers $m_{1,n}, m_{2,n},\cdots , m_{\left
\{ {n/2}\right \} ,n}$ given above, and for each $s=2,3,\cdots ,n-1$,
${\bf m}_{s}$ is a set of numbers $m_{1,s}, \cdots , m_{\left
\{ {s/2}\right \} ,s}$ and each $m_{rs}$, $r\le i$, runs independently the
values $h_{r,s}, h_{r,s}+1,\cdots , h_{r,s}+k-1$ and numbers
$m_{rs}$, $r> i$, run all half-integers
satisfying the same betweenness conditions as in
the case of irreducible representations of the nonclassical type
for $q$ not a root of unity (see [14]), that is, the conditions
$$
m_{i+1,2r+1}\ge m_{i+1,2r}\ge m_{i+2,2r+1}\ge m_{i+2,2r}\ge
\cdots \ge m_{r,2r+1}\ge m_{r,2r}\ge 1/2,
$$  $$
m_{i+1,2r}\ge m_{i+1,2r-1}\ge m_{i+2,2r}\ge
m_{i+2,2r-1}\ge \cdots \ge m_{r-1,2r-1}\ge m_{r,2r}\ge 1/2.
$$

To the set of numbers $\omega$ there corresponds the irreducible
finite dimensional representation $T''_\omega$ of the algebra
$U'_q({\rm so}_n)$.

If $p\le i$, then the operators $T''_\omega (I_{2p+1,2p})$ and
$T''_\omega (I_{2p,2p-1})$ act on the basis elements $| \xi \rangle$
by formulas (23) and (24) (if to replace in (24) $p$ by $p+1$).
If $p\ge i+1$, then the operators $T''_\omega (I_{2p+1,2p})$ and
$T''_\omega (I_{2p+2,2p+1})$ act as
$$
T''_\omega (I_{2p+1,2p})
| \xi_n\rangle =
\delta_{m_{p,2p},1/2}\, \frac{\epsilon _{2p+1}}{q^{1/2}-q^{-1/2}} D_{2p}
(\xi _n) | \xi_n\rangle + \qquad\qquad\qquad
$$
$$  \qquad\qquad\qquad
+\sum^p_{j=1} \frac{c_{j,2p} A^j_{2p}(\xi_n)}
{q^{l_{j,2p}}-q^{-l_{j,2p}} }
            \vert (\xi_n)^{+j}_{2p}\rangle -
\sum^p_{j=1} \frac{c^{-1}_{j,2p}A^j_{2p}((\xi_n)^{-j}_{2p})}
{q^{l_{j,2p}}-q^{-l_{j,2p}}}
|(\xi_n)^{-j}_{2p}\rangle    ,
$$
$$
T''_\omega (I_{2p+2,2p+1})\vert \xi_n\rangle=
\sum^{p}_{j=1} \frac{c_{j,2p+1}B^j_{2p+1}(\xi_n)}
{[2 l_{j,2p+1}-1][l_{j,2p+1}]_+}
\vert (\xi_n)^{+j}_{2p+1} \rangle -  \qquad\qquad\qquad
$$
$$\qquad\qquad\qquad
-\sum^{p}_{j=1}\frac {c^{-1}_{j,2p+1} B^j_{2p+1}((\xi_n)^{-j}_{2p+1})}
{[2 l_{j,2p+1}-1][l_{j,2p+1}-1]_+}
\vert (\xi_n)^{-j}_{2p+1}\rangle
+ \epsilon _{2p+2} {\hat C}_{2p+1}(\xi_n)
\vert \xi_n \rangle .
$$
In these formulas, $c_{sj}\equiv 1$ if $s>i$,
$(\xi_n)^{\pm j}_{k}$ means the tableau (21)
in which $j$-th component $m_{j,k}$ is replaced
by $m_{j,k}\pm 1.$ Matrix elements
$A^j_{2p}$ and $B^j_{2p+1}$ are given by the same formulas
as in (23) and (24), and
$$
{\hat C}_{2p+1}(\xi_n) = {
\prod_{s=1}^{p+1} [ l_{s,2p+2} ]_+
\prod_{s=1}^{p} [ l_{s,2p} ]_+ \over
\prod_{s=1}^{p} [l_{s,2p+1}]_+ [l_{s,2p+1} - 1]_+ } .
$$
$$
D_{2p} (\xi _n)=
\frac{\prod_{i=1}^p
[l_{i,2p+1}-\frac 12 ] \prod_{i=1}^{p-1} [l_{i,2p-1}-\frac 12 ] }
{\prod_{i=1}^{p-1}
[l_{i,2p}+\frac 12 ] [l_{i,2p}-\frac 12 ] } .
$$
In these formulas
$$
[b]_+=\frac{q^b+q^{-b}}{q-q^{-1}}.
$$

The fact that the operators $T''_\omega(I_{j,j-1})$ satisfy the defining
relations (1)--(3) of the algebra $U'_q({\rm so}_n)$ and irreducibility
of the representations $T''_\omega$ are proved in the same way as in
section 5.

\bigskip
{\sc Institute for Theoretical Physics, Ukrainian National Academy of
Sciences, Metrologichna str., 14-b, 03143, Kiev, Ukraine}

{\it E-mail address}: mmtpitp@bitp.kiev.ua

\bigskip
{\sc Institute for Theoretical Physics, Ukrainian National Academy of
Sciences, Metrologichna str., 14-b, 03143, Kiev, Ukraine}

{\it E-mail address}: aklimyk@gluk.apc.org

\end{document}